\documentclass[11pt]{article}
\usepackage{amsfonts}
\usepackage{color}
\newtheorem{theo}{Theorem}[section]
\newtheorem{lem}[theo]{Lemma}
\newtheorem{cor}[theo]{Corollary}

\newtheorem{defi}{Definition}[section]

\newcommand{\mysection}[1]{\section{#1} \setcounter{equation}{0}}
\newcommand{\proof}{{\sc Proof.} \quad}
\newcommand{\proofc}{{\sc Proof} \ }
\newcommand{\be}{\begin{equation} \label}
\newcommand{\ee}{\end{equation}}
\newcommand{\bea}{\begin{eqnarray}\label}
\newcommand{\eea}{\end{eqnarray}}
\newcommand{\bas}{\begin{eqnarray*}}
\newcommand{\eas}{\end{eqnarray*}}
\newcommand{\bit}{\begin{itemize}}
\newcommand{\eit}{\end{itemize}}
\newcommand{\qed}{\hfill$\Box$ \vskip.2cm}
\newcommand{\nn}{\nonumber}
\newcommand{\R}{\mathbb{R}}

\newcommand{\pO}{\partial\Omega}

\newcommand{\eps}{\varepsilon}

\newcommand{\hra}{\hookrightarrow}

\newcommand{\abs}{\\[5pt]}

\newcommand{\io}{\int_\Omega}
\newcommand{\bom}{\overline{\Omega}}

\newcommand{\dw}{D_w}
\newcommand{\dz}{D_z}
\newcommand{\ouz}{\overline{u}_0}
\newcommand{\es}{\eps_\star}
\newcommand{\ess}{\eps_{\star\star}}
\newcommand{\esss}{\eps_{\star\star\star}}
\newcommand{\oa}{\widehat{a}}
\newcommand{\ow}{\widehat{w}}
\newcommand{\oz}{\widehat{z}}
\begin{document}
\title{Asymptotic stability of spatial homogeneity
in a haptotxis model for oncolytic virotherapy}
\author{Youshan Tao\footnote{taoys@sjtu.edu.cn}\\
{\small School of Mathematical Sciences, Shanghai Jiao Tong University,}\\
{\small Shanghai 200240, P.R.~China}
\and
Michael Winkler\footnote{michael.winkler@math.uni-paderborn.de}\\
{\small Institut f\"ur Mathematik, Universit\"at Paderborn,}\\
{\small 33098 Paderborn, Germany} }
\date{}
\maketitle
\begin{abstract}
\noindent
  This work considers a model for oncolytic virotherapy, as given by the reaction-diffusion-taxis system
   \bas
        \left\{ \begin{array}{l}
        u_t = \Delta u - \nabla \cdot (u\nabla v)-\rho uz, \\[1mm]
        v_t = - (u+w)v, \\[1mm]
        w_t = \dw \Delta w - w + uz, \\[1mm]
    	z_t = \dz \Delta z - z - uz + \beta w,
         \end{array} \right.
  \eas
  in a smoothly bounded domain $\Omega\subset\R^2$, with parameters $\dw>0, \dz>0, \beta>0$ and $\rho\ge 0$.\abs
  Previous analysis has asserted that for all reasonably regular initial data, an associated no-flux type
  initial-boundary value problem admits a global classical solution, and that this solution is bounded if $\beta<1$,
  whereas whenever $\beta>1$ and $\frac{1}{|\Omega|}\io u(\cdot,0)>\frac{1}{\beta-1}$, infinite-time blow-up occurs at least
  in the particular case when $\rho=0$.\abs
  In order to provide an appropriate complement to this, the present work reveals that for any $\rho\ge 0$ and
  arbitrary $\beta>0$, at each prescribed
  level $\gamma\in (0,\frac{1}{(\beta-1)_+})$ one can identify an $L^\infty$-neighborhood of the homogeneous distribution
  $(u,v,w,z)\equiv (\gamma,0,0,0)$ within which all initial data lead to globally bounded solutions that stabilize toward
  the constant equilibrium $(u_\infty,0,0,0)$ with some $u_\infty>0$.\abs
\end{abstract}
\newpage
\mysection{Introduction}
Oncolytic virus particles are engineered for killing cancer cells, but they are little harmful to healthy cells.
The virions selectively adhere to the surface of cancer cells, then enter tumor cells via endocytosis,
enlarge their quantity through replication, and eventually cause the death of tumor cells. Upon lysis of a tumor cell,
a lot of new viruses are released, and they continue to infect adjacent tumor cells; the above process will be repeated
until all tumor cells are eradicated.
Accordingly, the visionary objective in this field is that
due to the considerable replication competence of viruses,
appropriately arranged treatments might provide efficient alternatives to conventional chemothapy,
with all its limitations linked to drug transport
(\cite{swabb}, \cite{jain}), and viral therapy has indeed already been
used in several clinical trails (\cite{bischoff}, \cite{coffey}, \cite{martuza}, \cite{russell}).\abs
Nevertheless, oncolytic efficacy of this novel therapy is also limited, not only by virus clearance due to various immune responses
(\cite{alemany}), but also by physical barriers such as interstitial fluid pressure and extracellular matrix (ECM) deposit
(\cite{vaha}, \cite{wong}). In order to figure out the role of the ECM in the spatio-temporal dynamics of virus spread within
a macro tissue including cancer cells and evaluate general effectiveness of oncolytic virotherapy, the authors in \cite{eftimie}
introduced a reaction-diffusion-taxis model that addresses the interaction between both uninfected and infected cancer cells, as well as ECM and
oncolytic virus particles. \abs
By neglecting any possible growth of uninfected tumor cells and ECM, in this study we consider a simplified version
of an originally more comprehensive model proposed in \cite{eftimie},
and will hence subsequently be concerned with the initial-boundary value problem
\be{0}
        \left\{ \begin{array}{lcll}
    	u_t &=& \Delta u - \nabla \cdot (u\nabla v) -\rho uz,
    	& x\in\Omega, \ t>0, \\[1mm]
	v_t &=& - (u+w)v,
    	& x\in\Omega, \ t>0, \\[1mm]
	w_t &=& \dw \Delta w - w + uz,
    	& x\in\Omega, \ t>0, \\[1mm]
	z_t &=& \dz \Delta z - z - uz + \beta w,
    	& x\in\Omega, \ t>0, \\[1mm]
    	& & \hspace*{-15mm}
    	(\nabla u -u\nabla v)\cdot \nu=\frac{\partial w}{\partial\nu} = \frac{\partial z}{\partial\nu}=0,
    	& x\in\pO, \ t>0, \\[1mm]
   	& & \hspace*{-15mm}
    	u(x,0)=u_0(x),
    	\quad v(x,0)=v_0(x),
    	\quad w(x,0)=w_0(x),
	\quad z(x,0)=z_0(x),
    	& x\in\Omega,
        \end{array} \right.
\ee
in a smoothly bounded domain $\Omega\subset\R^2$, with $\beta>0, \dw>0, \dz>0$ and $\rho\ge 0$,
and with the unknown variables $u, w, z$ and $v$ denoting the population densities
of uninfected cancer cells, infected cancer cells, virus particles and ECM, respectively.
Here as a crucial assumption from \cite{eftimie} that marks a substantial difference between (\ref{0}) and related more classical
reaction-diffusion models for virus dynamics (\cite{komarova}, \cite{pruess_mmnp}), we emphasize the hypothesis that
uninfected cancer cells can bias their motion upward ECM gradients due to attraction by some macromolecules trapped in ECM;
due to the fact that the ECM does not move, the resulting cross-diffusive migration is toward a non-diffusible quantity
and hence of haptotaxis type.\abs
It is quite precisely this latter circumstance that brings about considerable challenges for the mathematical analysis of (\ref{0}),
especially when focusing on questions related to qualitative solution behavior.
Indeed, previous studies concerned with related haptotaxis systems have mainly concentrated on establishing
mere solution theories (\cite{walker_webb}, \cite{tao_wang_sima2009},
\cite{taowin_JDE2014}, \cite{stinner_surulescu_winkler_SIMA},
\cite{zhigun_surulescu_uatay}, \cite{surulescu_win_CMS}),
with the only few exceptions available in the literature addressing rather specific settings
(\cite{taowin_NON}, \cite{pang_wang_m3as2018}, \cite{cao_zamp2016}),
\cite{taowin_sima2015}, \cite{litcanu_cmr} \cite{morales_tello},
\cite{win_jmpa2018}, \cite{hillen_painter_win_M3AS}).
After all, the crucial first equation in (\ref{0}) accounts for an essentially superlinear dampening mechanism
at least when $\rho$ is positive, and
some precedents have indeed revealed some significantly stabilizing effects of either precisely identical (\cite{rodriguez_win},
\cite{taowin_crime}) or related superlinear zero-order degradation terms in contexts of haptotaxis systems
(\cite{walker_webb}, \cite{taowin_JDE2014}, \cite{morales_tello},
\cite{stinner_surulescu_winkler_SIMA}, \cite{pang_wang_m3as2018});
in fact, it has recently been shown that if an additional logistic-type influence in the style of
an extra summand $\mu u(1-u)$ in the first equation from (\ref{0}), then solutions remain bounded and, if moreover $\beta<1$,
even approach the constant state $(1,0,0,0)$ in the large time limit (\cite{zhen_chen}).
In the absence of such further mechanisms, the identification of possible relaxing effects potentially induced by the
accordingly remaining and somewhat weaker dampening term $-\rho uz$ seems much less obvious,
especially in view of the coupling to the reaction-diffusion subsystem for $w$ and $z$, which at least in the case $\beta>1$
may itself apparently exhibit some strong tendency toward destibilization when forced by some appropriate $u$
(\cite{taowin_262}).\abs
Correspondingly, beyond a basic result on global smooth solvability
in widely arbitrary parameter settings and for all suitably regular initial data (\cite{taowin_261}),
the knowledge so far available for (\ref{0}) seems restricted to findings either addressing qualitative features
less subtle than precise convergence properties, e.g.~in the style of pointwise lower bounds or even unboundedness phenomena,
or concentrating on parameter constellations in which said $(w,z)$-subsystem remains subcritical in an appropriate sense.
Specifically, in \cite{taowin_265} it was seen that when $\rho=1$,
the size of $\beta$ relative to the value $\beta=1$ appears to be
decisive with regard to the question whether or not solutions may persistently remain above arbitrarily large levels
in their cancer cell population component $u$ throughout evolution, and that hence for appropriate efficiency of virotherapy
it might be advisable to assert virus reproduction rates fulfilling $\beta>1$.
A yet more drastic phenomenon indicating criticality of $\beta=1$ has been revealed in the borderline case $\rho=0$,
in which, namely, solutions to (\ref{0}) must be unbounded whenever $\beta>1$ and
$\frac{1}{|\Omega|} \io u_0 > \frac{1}{\beta-1}$, whereas in the semitrivial case when $v_0\equiv 0$, assuming that either
$\beta\le 1$, or $\frac{1}{|\Omega|} \io u_0 < \frac{1}{\beta-1}$, leads to globally bounded solutions (\cite{taowin_262}).
In extension of the latter, it has recently been found that also for general $\rho\ge 0$ and arbitrary $v_0$,
the corresponding solution of (\ref{0}) remains bounded when $\beta<1$ (\cite{taowin_263}).
Apart from this, some slightly more comprehensive variants of (\ref{0}) have been considered in \cite{yifu_wang} and \cite{xueyan_tao},
where in accordance with one of the models proposed in \cite{eftimie}
the inclusion of two further haptotaxis mechanisms, both of infected tumor cells and virions,
has been studied with respect to aspects of classical solvability and boundedness in the presence of certain suitably strong further
zero-order degradation (\cite{yifu_wang}), and of global smooth solvability in spatially one-dimensional settings
(\cite{xueyan_tao}).\abs
{\bf Main result.} \quad
According to the above, except for the case when $v\equiv 0$ and hence any tactic migration actually is absent,
in the case $\beta>1$ the dynamical features of (\ref{0}) seem widely unexplored in any planar domain;
in fact, the existing literature apparently even leaves open the question whether at all
some nontrivial bounded solutions can be found
in the presence of such supercritical virus production rates.
The purpose of the present study is to develop a method capable of asserting that this can indeed be achieved
for a considerably large set of initial data which are located in some neighborhood of certain spatially homogeneous distributions.
To make this more precise, importing the precise framework underlying the basic theory from \cite{taowin_261}
we shall henceforth assume that $\Omega\subset\R^2$ is a bounded domain with smooth boundary, and that
\be{init}
        \left\{ \begin{array}{l}
    	\mbox{$u_0, v_0, w_0$ and $z_0$ are nonnegative functions from $\bigcup_{\vartheta\in (0,1)} C^{2+\vartheta}(\bom)$,} \\[1mm]
    	\mbox{with $u_0\not \equiv 0$, $w_0\not\equiv 0$, $z_0\not\equiv 0$, $\sqrt{v_0} \in W^{1,2}(\Omega)$ and
	$\frac{\partial u_0}{\partial\nu}=\frac{\partial v_0}{\partial\nu}=
	\frac{\partial w_0}{\partial\nu}=\frac{\partial z_0}{\partial\nu}=0$ on $\pO$,}
    \end{array} \right.
\ee
recalling that then (\ref{0}) admits a unique global classical solution (cf. also Lemma \ref{lem_basic} below).\abs
Our main result now makes sure that whenever the corresponding initial deviation from certain constant values is suitably small,
these solutions will remain globally bounded, and even approach some constant profiles asymptotically;
in fact, we shall see that a statement of this flavor can also be established in the borderline case $\rho=0$ in which
no zero-order degrading influence acts on the resepctive first solution components, and that in this latter case
the associated domain of attraction even contains arbitrarily large $w_0$ and $z_0$.
As our method actually applies to any choice of $\beta>0$, the following formulation includes the range $\beta<1$
in which, in fact, a slightly more comprehensive result has already been achieved in \cite{taowin_263}:
\begin{theo}\label{theo23}
  Let $\Omega\subset\R^2$ be a bounded
  domain with smooth boundary, and let $\beta>0$ and $\gamma \in (0,\frac{1}{(\beta-1)_+})$.
  Then for each $M>0$ one can find
  $\eps=\eps(\beta,\gamma,M)>0$
  with the property that whenever $\rho\ge 0$ and $u_0, v_0, w_0$ and $z_0$ are such that (\ref{init}) holds and that
  \be{iu}
	\|u_0-\gamma\|_{L^\infty(\Omega)} <\eps
  \ee
  and
  \be{iv}
	\|v_0\|_{L^\infty(\Omega)} < \eps
  \ee
  as well as
  \be{iw}
	\|w_0\|_{L^\infty(\Omega)} < \min \Big\{ \frac{\eps}{\rho} \, , \, M \Big\},
  \ee
  and
  \be{iz}
	\|z_0\|_{L^\infty(\Omega)} < \min \Big\{ \frac{\eps}{\rho} \, , \, M \Big\},
  \ee
  there exists $u_\infty>0$ such that solution of (\ref{0}) satisfies
  \be{23.3}
	u(\cdot,t) \to u_\infty
	\qquad \mbox{in } L^\infty(\Omega)
  \ee
  and
  \be{23.4}
	v(\cdot,t) \to 0
	\qquad \mbox{in } L^\infty(\Omega)
  \ee
  as well as
  \be{23.5}
	w(\cdot,t) \to 0
	\qquad \mbox{in } L^\infty(\Omega)
  \ee
  and
  \be{23.6}
	z(\cdot,t) \to 0
	\qquad \mbox{in } L^\infty(\Omega)
  \ee
  as $t\to\infty$.
  Moreover, in the particular case when $\rho=0$ we have
  $u_\infty=\ouz$, so that for any $\gamma \in (0,\frac{1}{(\beta-1)_+})$
  the corresponding steady state solution $(\gamma,0,0,0)$
  of (\ref{0}) is asymptotically stable with respect to the norm
  in $(L^\infty(\Omega))^4$ in the above sense.
\end{theo}
{\bf Organization of the paper.} \quad
Our strategy will be based on the design of a self-map type reasoning, which presupposes a certain assumption on smallness
and decay of $z$ within an appropriate time interval, finally seen to actually be all of $[0,\infty)$ (Defintion \ref{dk}),
in order to derive appropriate boundedness and stabilization features of the solution as a whole,
which especially are consistent with said hypothesis.
Mainly due to fairly straightforward implications on pointwise lower bounds for $u$ and uniform decay of $v$
(Lemma \ref{lem25}), the core of our analysis will reveal that throughout the interval within which our assumption holds,
$w$ and $z$, along with a transformed version of $u$, form a subsolution to a cooperative parabolic system (Corollary \ref{cor26}).
Closing the loop of arguments will thus become possible through a derivation of upper estimates for $u$, $w$ and, particularly,
for $z$ by means of an associated comparison principle (Lemma \ref{lem28}).
On the basis of further temporally uniform regularity properties thereby implied, as seen through an appropriately organized bootstrap
procedure,		
thanks to a known conditional statement on stabilization in the first
solution component with respect to the norm in $L^2(\Omega)$ (Lemma \ref{lem52}) the outcome of this key step will be found  to
entail the claimed main result in Section \ref{sect5}.
\mysection{Preliminaries. Global existence and a conditional stabilization result for $u$}
%
%
%
%
%
%
%
%
Let us first recall that the conclusion of \cite{taowin_261} asserts global smooth solvability:
\begin{lem}\label{lem_basic}
  Let $\Omega\subset\R^2$ be a bounded
  domain with smooth boundary, let $\beta>0$ and $\rho\ge 0$, and suppose that
  $(u_0,v_0,w_0,z_0)$ satisfies (\ref{init}).
  Then the problem (\ref{0}) possesses a uniquely determined classical solution
  $(u,v,w,z) \in (C^{2,1}(\bom\times [0,\infty)))^4$ for which $v$ is nonnegative, and for which $u,w$ and $z$ are positive
  in $\bom\times (0,\infty)$.
\end{lem}
Following well-established basic strategy to conveniently reformulate the haptotactic interaction in (\ref{0}), as
widely used in related literature (\cite{fontelos_friedman_hu}, \cite{friedman_tello},
\cite{walker_webb}, \cite{taowin_JDE2014}), let us set
\be{a}
    a:= u e^{-v}
\ee
to see on the basis of (\ref{0}) that
\be{0a}
        \left\{ \begin{array}{lcll}
    	a_t &=& e^{-v} \nabla \cdot (e^v \nabla a)
	+ a(ae^v+w)v - \rho az,
    	& x\in\Omega, \ t>0, \\[1mm]
    	& & \hspace*{-13mm}
    	\frac{\partial a}{\partial\nu}=0,
    	& x\in\pO, \ t>0, \\[1mm]
    	& & \hspace*{-13mm}
    	a(x,0)=a_0(x):=u_0(x) e^{-v_0(x)},
    	& x\in\Omega.
        \end{array} \right.
\ee
In order to complete our small list of tokens imported from the literature, let us already here
recall from \cite{taowin_263} that appropriate assumptions on boundedness
of $(a,w,z)$ and decay of $v$ and $z$ are sufficient to ensure stabilization of $a$ in $L^2(\Omega)$.
This will be referred to in Section \ref{sect5} below.
\begin{lem}\label{lem52}
  Let $\beta>0$, $\gamma>0$ and $\rho\ge 0$, and suppose that for some initial data filfilling (\ref{init}),
  the corresponding solution of (\ref{0}) is such that
  \be{52.1}
	\sup_{t>0} \Big\{ \|a(\cdot,t)\|_{L^\infty(\Omega)}
      + \|w(\cdot,t)\|_{L^\infty(\Omega)}
	+ \|z(\cdot,t)\|_{L^\infty(\Omega)} \Big\} < \infty,
  \ee
  that
  \be{52.2}
  	v(\cdot, t)\to 0
  	\quad \mbox{in } L^\infty(\Omega)
  	\qquad\mbox{as } t\to \infty,
  \ee
  and that moreover
  \be{52.3}
	\int_0^\infty \io z <\infty.
  \ee
  Then there exists $u_\infty>0$ with the property that
  \bas
	\|u(\cdot,t)-u_\infty\|_{L^2(\Omega)} \to 0
	\qquad \mbox{as } t\to\infty.
  \eas
\end{lem}
\proof
  This can be seen by means of a verbatim copy of the reasoning from \cite[Lemma 6.5, Lemma 6.6 and Lemma 6.7]{taowin_263}.
\qed
\mysection{Some pointwise estimates for $u$ and $v$}
Next intending to set the frame for the announced self-map type argument, we first formulate an observation which,
although being quite elementary, contains the origin of our restriction on $\gamma$ in Theorem \ref{theo23}.
\begin{lem}\label{lem51}
  Let $\beta>0$ and $\gamma\in (0,\frac{1}{(\beta-1)_+})$.
  Then there exist $K_1(\beta,\gamma)>0$ and $\delta=\delta(\beta,\gamma)\in (0,1)$ such that
  \be{51.1}
	\frac{\gamma}{1-\delta} < K_1(\beta,\gamma) < \frac{\gamma+1-\delta}{\beta}.
  \ee
\end{lem}
\proof
  As our hypothesis $\gamma<\frac{1}{(\beta-1)_+}$ warrants that $\gamma<\frac{\gamma+1}{\beta}$, the number
  \bas
	K_1(\beta,\gamma):=\frac{1}{2} \Big( \gamma+ \frac{\gamma+1}{\beta}\Big)
  \eas
  satisfies $\gamma<K_1(\beta,\gamma) < \frac{\gamma+1}{\beta}$.
  Therefore, the claim follows by means of an argument based on continuous dependence.
\qed
We can thereby unambiguously formulate the core assumption underlying our subsequent analysis:
\begin{defi}\label{dk}
  Given $\beta>0$ and $\gamma\in (0,\frac{1}{(\beta-1)_+})$,
  we let
  \bas
	K_2(\beta,\gamma) := \max \Big\{ 1 \, , \, \frac{1}{K_1(\beta,\gamma)} \Big\} \ge 1,
  \eas
  where $K_1(\beta,\gamma)>0$ and $\delta(\beta,\gamma)\in (0,1)$ are as provided by Lemma \ref{lem51}.\abs
  Moreover, if $\rho\ge 0$ and $\eps>0$, and if $(u_0,v_0,w_0,z_0)$ are such that (\ref{init}) and (\ref{iu}) and (\ref{iv}) as well as
  \be{iw1}
	\rho \|w_0\|_{L^\infty(\Omega)} < \eps,
  \ee
  and
  \be{iz1}
	\rho \|z_0\|_{L^\infty(\Omega)} < \eps,
  \ee
  then we define
  \be{S}
	S(\beta,\gamma,\eps):= \bigg\{ T>0 \ \bigg| \
	\rho \|z(\cdot,t)\|_{L^\infty(\Omega)} < 2K_2(\beta,\gamma)\eps e^{-\delta(\beta,\gamma) t}
	\mbox{ for all } t\in (0,T) \bigg\}
  \ee
  and
  \be{T}
	T(\beta,\gamma,\eps):=\sup S(\beta,\gamma,\eps) \in (0,\infty].
  \ee
\end{defi}
A first and rather basic conclusion from the hypothesis included in (\ref{S}) can be obtained by a simple comparison argument.
\begin{lem}\label{lem21}
  Let $\beta>0$, $\gamma\in (0,\frac{1}{(\beta-1)_+})$ and $\rho\ge 0$, and assume that (\ref{init}) as well as (\ref{iu}), (\ref{iv}),
  (\ref{iw1}) and (\ref{iz1}) hold with some $\eps>0$.
  Then
  \be{21.1}
	u(x,t) \ge \Big\{ \min_{y\in\bom} a_0(y) \Big\} \cdot e^{-\frac{2K_2\eps}{\delta}}
	\qquad \mbox{for all $x\in\Omega$ and } t\in (0,T),
  \ee
  where $K_2=K_2(\beta,\gamma)>0$, $T=T(\beta,\gamma,\eps) \in (0,\infty]$ and
  $\delta=\delta(\beta,\gamma)\in (0,1)$ are as in Definition \ref{dk} and Lemma \ref{lem51}, respectively.
\end{lem}
\proof
  According to (\ref{0a}) and our definition of $S$,
  \bea{21.3}
	a_t &\ge& e^{-v} \nabla\cdot (e^v \nabla a) - \rho az \nn\\
	&\ge& e^{-v} \nabla\cdot (e^v \nabla a) - 2K_2\eps e^{-\delta t} a
	\qquad \mbox{in } \Omega\times (0,T),
  \eea
  and to derive a lower bound for $a$ from this, we let $c_1:=\min_{y\in\bom} a_0(y)$ and
  $\underline{a}(x,t):=\psi(t)$, $(x,t)\in\bom\times [0,\infty)$,
  where
  \bas
	\psi(t):= c_1 e^{ - \frac{2K_2\eps}{\delta} (1-e^{-\delta t})},
	\qquad t\ge 0.
%
  \eas
  As thus $\psi'(t)=-2K_2\eps e^{-\delta t} \psi(t)$ for $t>0$ and $\psi(0)=c_1$, it follows that
  \bas
	\underline{a}_t - e^{-v} \nabla \cdot (e^v \nabla \underline{a}) +2K_2\eps e^{-\delta t} \underline{a}
	= \psi'(t) + 2K_2\eps e^{-\delta t} \psi(t)=0
	\qquad \mbox{in } \Omega\times (0,T),
  \eas
  and that $\underline{a}(x,0) =c_1\le a_0(x)$ for all $x\in\Omega$, so that by means of the comparison principle
  we obtain that $\underline{a} \le a$ in $\Omega\times (0,T)$.
  Since $\psi(t) \ge c_1 e^{-\frac{2K_2\eps}{\delta}}$ for all $t\ge 0$, this entails (\ref{21.1}).
\qed
The following implication of the latter for the behavior of $v$ is quite obvious.
\begin{lem}\label{lem22}
  If $\beta>0$, $\gamma\in (0,\frac{1}{(\beta-1)_+})$ and $\rho\ge 0$, and if (\ref{init}) as well as (\ref{iu}), (\ref{iv}),
  (\ref{iw1}) and (\ref{iz1}) hold with some $\eps>0$, then
  \be{22.1}
	v(x,t) \le \|v_0\|_{L^\infty(\Omega)} \cdot
		\exp \bigg\{ - \Big\{ \min_{y\in\bom} a_0(y)\Big\} \cdot e^{-\frac{2K_2\eps}{\delta}} \cdot t \bigg\}
	\qquad \mbox{for all $x\in\Omega$ and } t\in (0,T),
  \ee
  where again $K_2=K_2(\beta,\gamma)>0$, $T=T(\beta,\gamma,\eps) \in (0,\infty]$ and
  $\delta=\delta(\beta,\gamma)\in (0,1)$ are as in Definition \ref{dk} and Lemma \ref{lem51}.
\end{lem}
\proof
  In view of the nonnegativity of $v$ and $w$, (\ref{0}) together with (\ref{21.1}) implies that
  $v_t= -(u+w)v \le -c_1 v$ in $\Omega\times (0,T)$, where
  $c_1:=\big\{ \min_{y\in\bom} a_0(y) \big\} \cdot e^{-\frac{2K_2\eps}{\delta}}$.
  Hence,
  \bas
	v(x,t) \le v_0(x) e^{-c_1 t} \le \|v_0\|_{L^\infty(\Omega)} e^{-c_1 t}
	\qquad \mbox{for all $x\in\Omega$ and } t\in (0,T),
  \eas
  which establishes (\ref{22.1}),  as claimed.
\qed
\mysection{Boundedness of $u$ and decay of $(w,z)$}
This section contains the core of our analysis by deriving and adequately exploiting the cooperative parabolic system (\ref{1a})
in order to establish suitable upper bounds for $u$ as well as for $w$ and $z$.
As a first step toward this, we refine the pointwise estimates from Lemma \ref{lem21} and Lemma \ref{lem22} by
now imposing an appropriate smallness assumption on the parameter $\eps$ in our hypotheses
(\ref{iu}), (\ref{iv}), (\ref{iw1}) and (\ref{iz1}).
\begin{lem}\label{lem25}
  Let $\beta>0$ and $\gamma\in (0,\frac{1}{(\beta-1)_+})$.
  Then there exists $\es(\beta,\gamma)>0$ such that if $\rho\ge 0$, and if
  (\ref{init}), (\ref{iu}) and (\ref{iv}) as well as
  (\ref{iw1}) and (\ref{iz1}) hold with some $\eps\in (0,\es(\beta,\gamma))$, then
  \be{25.1}
	u(x,t) \ge \gamma-\Big( \frac{2K_2}{\delta} \gamma+\gamma+2\Big) \cdot \eps
	\qquad \mbox{for all $x\in\Omega$ and } t\in (0,T)
  \ee
  and
  \be{25.2}
	v(x,t) \le \eps e^{-\frac{\gamma}{2}t}
	\qquad \mbox{for all $x\in\Omega$ and } t\in (0,T).
  \ee
  Here, as before, $K_2=K_2(\beta,\gamma)>0$, $T=T(\beta,\gamma,\eps) \in (0,\infty]$ and
  $\delta=\delta(\beta,\gamma)\in (0,1)$ are taken from Definition \ref{dk} and Lemma \ref{lem51}.
\end{lem}
\proof
  By l'Hospital's rule, abbreviating $c_1:=\frac{2K_2}{\delta}+1$, we see that
  \bas
	\frac{(c_1\gamma+2-e^{-c_1\eps}) \cdot\eps}{1-e^{-c_1 \eps}}
	\to \lim_{\eps\searrow 0} \frac{c_1\gamma+2-e^{-c_1\eps} + c_1\eps e^{-c_1 \eps}}{c_1 e^{-c_1\eps}}
	= \gamma+\frac{1}{c_1}>\gamma,
  \eas
  so that it is possible to pick $\eps_1=\eps_1(\beta,\gamma)>0$ in such a way that
  \bas
	(c_1\gamma+2- e^{-c_1\eps}) \cdot \eps
	\ge \gamma \cdot (1-e^{-c_1 \eps})
	\qquad \mbox{for all } \eps\in (0,\eps_1)
  \eas
  and hence
  \be{25.22}
	(\gamma-\eps) e^{-c_1 \eps} \ge \gamma-(c_1\gamma+2) \cdot \eps
	\qquad \mbox{for all } \eps\in (0,\eps_1).
  \ee
  Furthermore, observing that
  \bas
	(\gamma-\eps) e^{-c_1\eps} \to \gamma>\frac{\gamma}{2}
	\qquad \mbox{as } \eps\searrow 0,
  \eas
  we can fix $\eps_2=\eps_2(\beta,\gamma)>0$ fulfilling
  \be{25.23}
	(\gamma-\eps)e^{-c_1 \eps} \ge \frac{\gamma}{2}
	\qquad \mbox{for all } \eps\in (0,\eps_2).
  \ee
  Therefore, if we let
  \bas
	\es\equiv \es(\beta,\gamma):=\min \Big\{ \gamma, \eps_1 (\beta,\gamma) \, , \, \eps_2(\beta,\gamma) \Big\}
  \eas
  and assume (\ref{init}), (\ref{iu}), (\ref{iv}), (\ref{iw1}) and (\ref{iz1}) with some $\eps\in (0,\es)$, then since
  \be{25.3}
	a_0(x) = u_0(x) e^{-v_0(x)} \ge (\gamma-\eps) e^{-\eps}
	\qquad \mbox{for all } x\in\Omega
  \ee
  by (\ref{iu}), (\ref{iv}) and (\ref{a}), and since thus
  \bas
	u(x,t) \ge (\gamma-\eps) e^{-\eps} e^{-\frac{2K_2\eps}{\delta}} = (\gamma-\eps) e^{-c_1 \eps}
	\qquad \mbox{for all $x\in\Omega$ and } t\in (0,T)
  \eas
  by Lemma \ref{lem21} and the fact that $\es\le \gamma$,
  using the restriction $\es\le \eps_1$ we firstly obtain (\ref{25.1}) as a consequence of (\ref{25.22}).
  Apart from that, the inequality $\es\le \eps_2$ in conjunction with Lemma \ref{lem22} and (\ref{iv}) ensures that
  (\ref{25.3}), secondly, guarantees that
  \bas
	v(x,t) &\le& \eps \exp \Big\{ -(\gamma-\eps) e^{-\eps} \cdot e^{-\frac{2K_2\eps}{\delta}} \cdot t\Big\} \\
	&=& \eps \exp \Big\{ -(\gamma-\eps) e^{-c_1 \eps} \cdot t\Big\} \\
	&\le& \eps e^{-\frac{\gamma}{2} t}
	\qquad \mbox{for all $x\in\Omega$ and } t\in (0,T)
  \eas
  thanks to (\ref{25.23}).
\qed
For such choices of $\eps$ this directly entails that, indeed, the triple $(a,w,z)$ forms a subsolution of a cooperative
reaction-diffusion system.
\begin{cor}\label{cor26}
  Let $\beta>0$ and $\gamma\in (0,\frac{1}{(\beta-1)_+})$, and let $\es(\beta,\gamma)>0$ be as in Lemma \ref{lem25}.
  Then whenever $\rho\ge 0$ and
  (\ref{init}), (\ref{iu}) and (\ref{iv}) as well as
  (\ref{iw1}) and (\ref{iz1}) hold with some $\eps\in (0,\es(\beta,\gamma))$, the solution of (\ref{0}) satisfies
  \be{1a}
        \left\{ \begin{array}{lcll}
    	a_t &\le& e^{-v} \nabla \cdot (e^v \nabla a)
	+ \eps e^\eps e^{-\frac{\gamma}{2}t} a^2 + \eps e^{-\frac{\gamma}{2}t} aw,
    	& x\in\Omega, \ t>0, \\[1mm]
    	w_t &\le& \dw \Delta w - w + e^\eps az,
    	& x\in\Omega, \ t>0, \\[1mm]
	z_t &\le& \dz \Delta z - \Big\{ \gamma +1 - \Big(\frac{2K_2}{\delta} \gamma+\gamma+2\Big) \cdot \eps \Big\} \cdot z + \beta w,
    	& x\in\Omega, \ t\in (0,T),
        \end{array} \right.
  \ee
  where once more $K_2=K_2(\beta,\gamma)>0$, $T=T(\beta,\gamma,\eps) \in (0,\infty]$ and
  $\delta=\delta(\beta,\gamma)\in (0,1)$ are given by Definition \ref{dk} and Lemma \ref{lem51}.
\end{cor}
\proof
  In (\ref{0a}), we only need to use (\ref{a}), (\ref{25.1}) and (\ref{25.2}) to see that thanks to the fact that $e^v \le e^\eps$
  in $\Omega\times (0,T)$ by the latter, we have
  \bas
	 (ae^v +w) v
	\le a(ae^\eps +w) \cdot \eps e^{-\frac{\gamma}{2}t}
  \eas
  and
  \bas
	-w+uz = -w+ae^v z \le -w + e^\eps az
  \eas
  as well as
  \bas
	-z-uz + \beta w \le -z - \bigg\{ \gamma - \Big(\frac{2K_2}{\delta} \gamma+\gamma+2\Big) \cdot \eps \bigg\} \cdot z +\beta w
  \eas
  in $\Omega\times (0,T)$.
\qed
Our construction of an appropriate supersolution to (\ref{1a}) will, in its crucial first component, involve
a spatially homogeneous time-dependent function taken from the family of solutions to quadratically forced Bernoulli-type
ODE problems addressed in the following lemma (cf.~(\ref{oa}) below).
\begin{lem}\label{lem27}
  Let $\gamma>0$ and $A>0$.
  Then there exists $\ess=\ess(\gamma,A)>0$ such that if $\eps\in (0,\ess)$, then the initial value problem
  \be{0phi}
	\left\{ \begin{array}{l}
	\varphi'(t) = \eps e^\eps e^{-\frac{\gamma}{2}t} \varphi^2(t)
	+ A\eps e^{-\frac{\gamma}{2}t} \varphi(t),
	\qquad t>0, \\[2mm]
	\varphi(0)=\gamma+\eps,
	\end{array} \right.
  \ee
  possesses a globally defined solution $\varphi\in C^1([0,\infty))$ fulfilling
  \be{27.1}
	\varphi(t) \le \gamma + (16\gamma+72A+1) \cdot \eps
	\qquad \mbox{for all } t>0.
  \ee
\end{lem}
\proof
  Given $\gamma>0$ and $A>0$, we abbreviate $c_1:=16\gamma+72A+1$ and let
  \be{27.2}
	\ess=\ess(\gamma,A):=\min \Big\{ \ln 2 \, , \, \frac{\gamma}{c_1} \, , \, \frac{\gamma \ln 2}{2A} \Big\}.
  \ee
  Then assuming that $\eps\in (0,\ess)$ and letting $\varphi\in C^1([0,T_\eps))$ denote tha corresponding solution of (\ref{0phi}),
  extended up to its maximal existence time $T_\eps \in (0,\infty]$, by explicitly solving (\ref{0phi}) we see that
  \be{27.3}
	\frac{1}{\varphi(t)}
	= \frac{1}{\gamma+\eps} \cdot \exp \bigg\{ -A\eps \int_0^t e^{-\frac{\gamma}{2}s} ds \bigg\}
	- \eps e^\eps \int_0^t \exp \bigg\{ -A\eps \int_s^t e^{-\frac{\gamma}{2}\sigma} d\sigma \bigg\} \cdot e^{-\frac{\gamma}{2}s} ds
	\qquad \mbox{for all } t\in (0,T_\eps).
  \ee
  Here since
  \be{27.4}
	\int_0^t e^{-\frac{\gamma}{2}s} ds \le \frac{2}{\gamma}
	\qquad \mbox{for all } t>0,
  \ee
  we obtain that
  \bas
	\frac{1}{\gamma+\eps} \cdot \exp \bigg\{ -A\eps \int_0^t e^{-\frac{\gamma}{2}s} ds \bigg\}
	\ge \frac{1}{\gamma+\eps} \cdot e^{-\frac{2A\eps}{\gamma}}
	\qquad \mbox{for all } t>0,
  \eas
  whereas simply estimating
  \bas
	\int_s^t e^{-\frac{\gamma}{2}\sigma} d\sigma \ge 0
	\qquad \mbox{for all $s\ge 0$ and } t\ge s,
  \eas
  noting that $e^\eps \le 2$ by (\ref{27.2}) we find that again due to (\ref{27.4}),
  \bas
	\eps e^\eps \int_0^t \exp \bigg\{ -A\eps \int_s^t e^{-\frac{\gamma}{2}\sigma} d\sigma \bigg\} \cdot e^{-\frac{\gamma}{2}s} ds
	&\le& 2\eps \int_0^t e^{-\frac{\gamma}{2}s} ds \\
	&\le& \frac{4\eps}{\gamma}
	\qquad \mbox{for all } t>0.
  \eas
  Therefore, (\ref{27.3}) entails that
  \bea{27.5}
	e^\frac{2A\eps}{\gamma} \cdot \Big\{ \frac{\gamma+c_1\eps}{\varphi(t)} -1 \Big\}
	&\ge& e^\frac{2A\eps}{\gamma} \cdot \bigg\{
		(\gamma+c_1\eps) \cdot \Big\{ \frac{1}{\gamma+\eps} \cdot e^{-\frac{2A\eps}{\gamma}} - \frac{4\eps}{\gamma} \Big\} -1
		\bigg\} \nn\\
	&=& \frac{\gamma+c_1\eps}{\gamma+\eps}
	- e^\frac{2A\eps}{\gamma} \cdot \Big\{ 1 + \frac{4(\gamma+c_1\eps)}{\gamma} \cdot \eps \Big\}
	\qquad \mbox{for all } t\in (0,T_\eps),
  \eea
  where since $\eps\le \gamma$ by (\ref{27.2}),
  \be{27.6}
	\frac{\gamma+c_1\eps}{\gamma+\eps}
	= 1+ \frac{c_1-1}{\gamma+\eps} \cdot\eps
	\ge 1+ \frac{c_1-1}{2\gamma} \cdot \eps,
  \ee
  because $c_1>1$.
  Furthermore, relying on the fact that $e^s \le 1+2s$ for all $s\in [0,\ln 2]$ we can make use of the rightmost restriction contained
  in (\ref{27.2}) to see that since clearly $c_1\eps \le \gamma$ and $\eps\le 1$ by (\ref{27.2}),
  \bas
	e^\frac{2A\eps}{\gamma} \cdot \Big\{ 1 + \frac{4(\gamma+c_1\eps)}{\gamma} \cdot \eps \Big\}
	&\le& \Big( 1+\frac{4A\eps}{\gamma}\Big) \cdot (1+8\eps) \\
	&=& 1+\Big(\frac{4A}{\gamma}+8\Big) \cdot \eps + \frac{32A\eps}{\gamma} \cdot\eps \\
	&\le& 1+\Big( \frac{36A}{\gamma} +8\Big)\cdot \eps.
  \eas
  As thus, by (\ref{27.6}),
  \bas
	\frac{\gamma+c_1\eps}{\gamma+\eps}
	- e^\frac{2A\eps}{\gamma} \cdot \Big\{ 1 + \frac{4(\gamma+c_1\eps)}{\gamma} \cdot \eps \Big\}
	\ge \Big\{ \frac{c_1-1}{2\gamma} -\Big(\frac{36A}{\gamma}+8\Big) \Big\} \cdot\eps
	=0
  \eas
  according to our choice of $c_1$, from (\ref{27.5}) we infer that indeed $\varphi(t) \le \gamma+c_1\eps$ for all
  $t\in (0,T_\eps)$, which implies that in fact we must have $T_\eps=\infty$, and that (\ref{27.1}) holds.
\qed
Based on the latter, we can now quite easily find spatially constant supersolutions to (\ref{1a}) with
properties favorable for our purposes, and thereby accomplish the main step in our reasoning.
\begin{lem}\label{lem28}
  Let $\beta>0$, $\gamma\in (0,\frac{1}{(\beta-1)_+})$ and $M>0$.
  Then there exist $\esss=\esss(\beta,\gamma,M)>0$ and $C=C(\beta,\gamma,M)>0$ with the property that if $\rho\ge 0$
  and (\ref{init}) as well as (\ref{iu})-(\ref{iz}) hold with some $\eps\in (0,\esss)$, we have
  \be{28.1}
	u(x,t) \le C	
	\qquad \mbox{for all $x\in\Omega$ and } t \in (0,T)
  \ee
  and
  \be{28.2}
	w(x,t) \le K_1 K_2 \cdot \min\Big\{ \frac{\eps}{\rho} \, , \, M \Big\} \cdot e^{-\delta t}
	\qquad \mbox{for all $x\in\Omega$ and } t \in (0,T)
  \ee
  as well as
  \be{28.3}
	z(x,t) \le K_2 \cdot \min\Big\{ \frac{\eps}{\rho} \, , \, M \Big\} \cdot e^{-\delta t}
	\qquad \mbox{for all $x\in\Omega$ and } t \in (0,T),
  \ee
  with $K_i=K_i(\beta,\gamma)>0$, $i\in \{1,2\}$, $T=T(\beta,\gamma,\eps) \in (0,\infty]$ and
  $\delta=\delta(\beta,\gamma)\in (0,1)$ taken from Definition \ref{dk} and Lemma \ref{lem51}.
\end{lem}
\proof
  As Lemma \ref{lem51} warrants that $\frac{\gamma}{1-\delta} < K_1 <\frac{\gamma+1-\delta}{\beta}$,
  given $M>0$ we can fix $\eps_1=\eps_1(\beta,\gamma,M) \in (0,1)$ such that
  \be{28.11}
	\frac{(\gamma+c_1 \eps_1) e^{\eps_1}}{1-\delta} \le K_1 \le
	\frac{\gamma+1 - c_2 \eps_1 - \delta}{\beta},
  \ee
  where we have set
  \be{28.111}
	c_1\equiv c_1(\beta,\gamma,M):=16\gamma+72K_1 K_2 M + 1
	\qquad \mbox{and} \qquad
	c_2\equiv c_2(\beta,\gamma):=\frac{2K_2(\beta,\gamma)}{\delta(\beta,\gamma)} \cdot \gamma + \gamma+2.
  \ee
  Then letting
  \be{28.12}
	\esss\equiv \esss(\beta,\gamma,M):=
	\min \Big\{ \eps_1(\beta,\gamma,M) \, , \, \es(\beta,\gamma) \, , \, \ess(\gamma, K_1 K_2 M) \Big\}
  \ee
  with $\es(\cdot,\cdot)>0$ and $\ess(\cdot,\cdot)$ taken from Lemma \ref{lem25} and Lemma \ref{lem27},
  we henceforth assume that $\rho\ge 0$, and that (\ref{init}) and (\ref{iu})-(\ref{iz}) are satisfied with some $\eps\in (0,\esss)$.
  We then abbreviate
  \be{28.4}
	B:=\max \Big\{ \|z_0\|_{L^\infty(\Omega)} \, , \, \frac{1}{K_1} \cdot \|w_0\|_{L^\infty(\Omega)} \Big\}
  \ee
  as well as
  \be{28.5}
	A:=K_1 B,
  \ee
  and first note that then due to (\ref{iw}), (\ref{iz}) and our definition of $K_2$,
  \be{28.55}
	B < K_2 M
	\qquad \mbox{and} \qquad
	A < K_1 K_2 M
  \ee
  as well as
  \be{28.56}
	B < \frac{K_2\eps}{\rho}
	\qquad \mbox{and} \qquad
	A < \frac{K_1 K_2 \eps}{\rho}.
  \ee
  In order to construct a corresponding supersolution triple $(\oa,\ow,\oz)$, based on the above choices we let
  \be{oa}
	\oa(x,t):=\varphi(t),
	\qquad x\in\bom, \ t\ge 0,
  \ee
  as well as
  \bas
	\ow(x,t):=Ae^{-\delta t},
	\qquad x\in\bom, \ t\ge 0,
  \eas
  and
  \bas
	\oz(x,t):=Be^{-\delta t},
	\qquad x\in\bom, \ t\ge 0,
  \eas
  where in accordance with Lemma \ref{lem27}, $\varphi\in C^1([0,\infty))$ denotes the solution of (\ref{0phi}). Then thanks to the
  latter,
  \bea{28.7}
	& & \hspace{-20mm}
	\oa_t - e^{-v} \nabla \cdot (e^v \nabla\oa)
	- \eps e^\eps e^{-\frac{\gamma}{2}t} \oa^2
	- \eps e^{-\frac{\gamma}{2}t} \oa\ow \nn\\
	&=& \varphi'(t) - \eps e^\eps e^{-\frac{\gamma}{2}t} \varphi^2(t)
	-A\eps e^{-\frac{\gamma}{2}t} e^{-\delta t} \varphi(t) \nn\\
	&\ge & \varphi'(t) - \eps e^\eps e^{-\frac{\gamma}{2}t} \varphi^2(t)
	-A\eps e^{-\frac{\gamma}{2}t} \varphi(t) \nn\\[2mm]
	&=& 0
	\qquad \mbox{in } \Omega\times (0,T),
  \eea
  while the upper estimate in (\ref{27.1}) along with the second inequality in (\ref{28.55}) and our definition of $c_1$ guarantees that
  \bea{28.8}
	\ow_t - \dw \Delta \ow + \ow - e^\eps \oa \oz
	&=& - \delta A e^{-\delta t} + A e^{-\delta t} - e^\eps B e^{-\delta t} \varphi \nn\\
	&\ge & - \delta A e^{-\delta t} + A e^{-\delta t} - (\gamma+c_1\eps) e^\eps B e^{-\delta t} \nn\\
	&=& \Big\{ (1-\delta)A - (\gamma+c_1\eps) e^\eps B \Big\} \cdot e^{-\delta t} \nn\\[2mm]
	&\ge& 0
	\qquad \mbox{in } \Omega\times (0,\infty),
  \eea
  because by (\ref{28.5}) and the left inequality in (\ref{28.11}),
  \bas
	(1-\delta)A - (\gamma+c_2\eps) e^\eps B
	&=& \Big\{ (1-\delta)c_1 - (\gamma+c_2\eps) e^\eps \Big\} \cdot B \\
	&\ge& \Big\{ (1-\delta)c_1 - (\gamma+c_2\eps_1) e^{\eps_1} \Big\} \cdot B \\[2mm]
	&\ge& 0.
  \eas
  Likewise,
  \bea{28.9}
	& & \hspace*{-35mm}
	\oz_t - \dz \Delta\oz + \Big\{ \gamma +1 - \Big(\frac{2K_2}{\delta} \gamma+\gamma+2\Big) \cdot \eps \Big\} \cdot \oz
	- \beta \ow \nn\\
	&=& -\delta B e^{-\delta t}
	+ (\gamma +1 - c_2 \eps) \cdot B e^{-\delta t} - \beta A e^{-\delta t}
		\nn\\
	&=& (\gamma+1 - c_2\eps - \delta - K_1\beta) Be^{-\delta t} \nn\\[2mm]
	&\ge& 0
	\qquad \mbox{in } \Omega\times (0,\infty),
  \eea
  for
  \bas
	\gamma+1-c_2\eps -\delta-K_1 \beta
	\ge \gamma+1-c_2\eps_1 -\delta-K_1 \beta \ge 0
  \eas
  due to the right inequality in (\ref{28.11}). Since, apart from that,
  \bas
	a_0(x) = u_0(x) e^{-v_0(x)}
	\le u_0(x)
	\le \gamma+\eps
	= \varphi(0)
	= \oa(x,0)
	\qquad \mbox{for all } x\in\Omega
  \eas
  by (\ref{iu}), and since (\ref{28.4}) and (\ref{28.5}) guarantee that
  \bas
	w_0(x) \le \|w_0\|_{L^\infty(\Omega)}
	\le K_1 B = A = \ow(x,0)
	\qquad \mbox{for all } x\in\Omega
  \eas
  and
  \bas
	z_0(x) \le \|z_0\|_{L^\infty(\Omega)} \le B = \oz(x,0)
	\qquad \mbox{for all } x\in\Omega,
  \eas
  we may make use of the fact that the parabolic system in (\ref{1a}) is cooperative to conclude from a corresponding comparison
  principle that
  \bas
	a\le \oa, \quad
	w\le \ow
	\quad \mbox{and} \quad
	z\le \oz
	\qquad \mbox{in } \Omega\times (0,\infty).
  \eas
  By definition of $\ow$ and $\oz$, in view of (\ref{28.55}) and (\ref{28.56})
  the two latter inequalities directly yield (\ref{28.2}) and (\ref{28.3}), while (\ref{28.1})
  is a consequence of the bounds for $\varphi$ and $v$ asserted by Lemma \ref{lem27} and Lemma \ref{lem25}.
\qed
Through (\ref{28.3}), the latter especially enables us to close the loop implictly opened in Definition \ref{dk}.
\begin{cor}\label{cor288}
  Let $\beta>0$, $\gamma\in (0,\frac{1}{(\beta-1)_+})$, $\rho\ge 0$ and $M>0$, and suppose
  that (\ref{init}) and (\ref{iu})-(\ref{iz}) hold with some $\eps\in (0,\esss)$, where $\esss=\esss(\beta,\gamma,M)>0$
  is as given by Lemma \ref{lem28}. Then in Definition \ref{dk} we have $T(\beta,\gamma,\eps)=\infty$.
\end{cor}
\proof
  Since (\ref{28.3}) particularly entails that
  \bas
	\rho \|z(\cdot,t)\|_{L^\infty(\Omega)} \le K_2(\beta,\gamma) \eps e^{-\delta(\beta,\gamma) t}
	\qquad \mbox{for all } t\in (0,T(\beta,\gamma,\eps)),
  \eas
  assuming that $T(\beta,\gamma,\eps)$ be finite would readily lead to a contradiction to the defintion of $S(\beta,\gamma,\eps)$
  and the continuity of $z$.
\qed
\mysection{A global H\"older bound for $u$. Proof of Theorem \ref{theo23}}\label{sect5}
In view of Lemma \ref{lem28} and Corollary \ref{cor288}, it remains to be shown that the $L^2$ stabilization process in the first
solution component, as thus clearly asserted by Lemma \ref{lem52}, in fact can be turned into the uniform convergence statement
in (\ref{23.3}).
This will be achieved on the basis on the following boundedness property of the haptotactic gradient in (\ref{0}),
as resulting from a series of testing procedures applied to the first three equations therein.
\begin{lem}\label{lem33}
  Let $\beta>0$, $\gamma\in (0,\frac{1}{(\beta-1)_+})$, $\rho\ge 0$ and $M>0$, and assume
  (\ref{init}) and (\ref{iu})-(\ref{iz}) with some $\eps\in (0,\esss)$, and with $\esss=\esss(\beta,\gamma,M)>0$
  taken from Lemma \ref{lem28}.
  Then there exists $C>0$ such that
  \be{33.2}
	\io |\nabla v(\cdot,t)|^4  \le C
	\qquad \mbox{for all } t>0.
  \ee
\end{lem}
\proof
  Relying on Lemma \ref{lem28}, let us pick $c_1>0, c_2>0$ and $c_3>0$ such that
  \be{33.22}
	u(x,t) \le c_1,
	\quad
	w(x,t) \le c_2
	\quad \mbox{and} \quad
	z(x,t) \le c_3
	\qquad \mbox{for all $x\in\Omega$ and } t>0,
  \ee
  which due to (\ref{a}) and Lemma \ref{lem22} particularly means that
  \bas
	f=f(x,t):=a(ae^v+w)v,
	\qquad x\in\Omega \, t>0,
  \eas
  satisfies
  \bas
	|f(x,t)| \le c_5:=c_1 \cdot (c_1+c_2) \cdot c_4
	\qquad \mbox{for all $x\in\Omega$ and } t>0
  \eas
  with $c_4:=\|v_0\|_{L^\infty(\Omega)}$.
  Therefore, testing the identity
  \bas
	a_t = \Delta a + \nabla v\cdot\nabla a + f(x,t),
	\qquad x\in\Omega, \ t>0,
  \eas
  as contained in (\ref{0a}), by $-\Delta a$ and using Young's inequality shows that
  \bea{33.4}
	\frac{1}{2} \frac{d}{dt} \io |\nabla a|^2 + \io |\Delta a|^2
	&=& - \io (\nabla v\cdot\nabla a)\Delta a
	- \io f\Delta a \nn\\
	&\le& \frac{1}{2} \io |\Delta a|^2
	+ \io |\nabla v\cdot\nabla a|^2 +  \io f^2 \nn\\
	&\le& \frac{1}{2} \io |\Delta a|^2
	+ \io |\nabla v|^2 |\nabla a|^2
	+ c_5^2 |\Omega|
	\qquad \mbox{for all } t>0.
  \eea
  Here we combine a Gagliardo-Nirenberg type interpolation with standard elliptic regularity theory to find $c_6>0$ fulfilling
  \be{33.5}
	\|\nabla\varphi\|_{L^4(\Omega)}^4 \le c_6\|\Delta\varphi\|_{L^2(\Omega)}^2 \|\varphi\|_{L^\infty(\Omega)}^2
	\qquad \mbox{for all $\varphi\in W^{2,2}(\Omega)$ such that $\frac{\partial\varphi}{\partial\nu}=0$ on } \pO,
  \ee
  which, again thanks to (\ref{33.22}) and (\ref{a}), firstly implies that
  \bas
	\frac{1}{2} \io |\Delta a|^2
	\ge \frac{1}{2c_1^2 c_6} \io |\nabla a|^4
	\qquad \mbox{for all } t>0.
  \eas
  Therefore, through two applications of Young's inequality we infer from (\ref{33.5}) that
  \bas
	\frac{1}{2} \frac{d}{dt} \io |\nabla a|^2
	+ \frac{1}{2c_1^2 c_6} \io |\nabla a|^4
	+ \io |\nabla a|^2
	&\le& \io |\nabla v|^2 |\nabla a|^2
	+ c_5^2 |\Omega| + \io |\nabla a|^2 \\
	&\le& \frac{1}{8c_1^2 c_6} \io |\nabla a|^4
	+ 2c_1^2 c_6 \io |\nabla v|^4
	+ c_5^2 |\Omega| \\
	& & + \frac{1}{8c_1^2 c_6} \io |\nabla a|^4
	+ 2c_1^2 c_6 |\Omega|
	\qquad \mbox{for all } t>0
  \eas
  and hence
  \be{33.7}
	\frac{d}{dt} \io |\nabla a|^2 + c_7 \io |\nabla a|^4 + 2\io |\nabla a|^2
	\le c_8 \io |\nabla v|^4 + c_8
	\qquad \mbox{for all } t>0
  \ee
  with $c_7:=\frac{1}{2c_1^2 c_6}$ and
  $c_8:=\max \big\{ 4c_1^2 c_6 \, , \, 2c_5^2 |\Omega| + 4c_1^2 c_6|\Omega|\big\}$.\abs
  In order to appropriately compensate the first summand on the right of (\ref{33.7}), we next use the second equation in
  (\ref{0a}) to see that
  \bea{33.8}
	\frac{1}{4} \frac{d}{dt} \io |\nabla v|^4
	&=& - \io |\nabla v|^2 \nabla v\cdot \nabla (ave^v + vw) \nn\\
	&=& - \io a(v+1) e^v |\nabla v|^4
	- \io ve^v |\nabla v|^2 \nabla v\cdot\nabla a \nn\\
	& & - \io w|\nabla v|^4
	- \io v|\nabla v|^2 \nabla v\cdot\nabla w
	\qquad \mbox{for all } t>0,
  \eea
  where the second last summand is nonpositive, where we moreover
  recall the uniform positivity statement for $u=ae^v$ from Lemma \ref{lem21}
  to pick $c_9 \in (0,4]$ fulfilling
  \be{33.86}
	\io a(v+1) e^v |\nabla v|^4 \ge \io ae^v |\nabla v|^4 \ge c_9 \io |\nabla v|^4
	\qquad \mbox{for all } t>1,
  \ee
  and where by Lemma \ref{lem22} and Young's inequality,
  \bea{33.87}
	- \io v|\nabla v|^2 \nabla v\cdot\nabla w
	&\le& c_4 \io |\nabla v|^3 |\nabla w| \nn\\
	&=& \io \Big\{ \frac{c_9}{2} |\nabla v|^4 \Big\}^\frac{3}{4} \cdot \Big\{ \Big(\frac{2}{c_9}\Big)^\frac{3}{4} c_4|\nabla w|\Big\}
		\nn\\
	&\le& \frac{c_9}{2} \io |\nabla v|^4
	+ \frac{8 c_4^4}{c_9^3} \io |\nabla w|^4
	\qquad \mbox{for all } t>0.
  \eea
  Now in estimating the second summand on the right-hand side of (\ref{33.8}) we proceed slightly more carefully in order to retain
  a potentially small factor: Indeed, given any $t_0\ge 0$ we may combine our definition of $c_4$ with Lemma \ref{lem22} to obtain that,
  again by Young's inequality,
  \bea{33.877}
	- \io v e^v |\nabla v|^2 \nabla v\cdot\nabla a
	&\le& \|v(\cdot,t_0)\|_{L^\infty(\Omega)} e^{c_4} \io |\nabla v|^3 |\nabla a| \nn\\
	&=& \io \Big\{ \frac{c_9}{4}|\nabla v|^4 \Big\}^\frac{3}{4} \cdot
	\Big\{ \Big(\frac{4}{c_9}\Big)^\frac{3}{4} e^{c_4} \|v(\cdot,t_0)\|_{L^\infty(\Omega)} |\nabla a| \Big\} \nn\\
	&\le& \frac{c_9}{4} \io |\nabla v|^4
	+ \frac{64 e^{4c_4}}{c_9^3} \|v(\cdot,t_0)\|_{L^\infty(\Omega)}^4 \io |\nabla a|^4
	\qquad \mbox{for all } t>t_0.
  \eea
  If we write $c_{10}:=\max\Big\{ \frac{32c_4^4}{c_9^3} \, , \, \frac{256 e^{4c_4}}{c_9^3}\Big\}$, from (\ref{33.8})-(\ref{33.877})
  we thus infer that whenever $t_0\ge 1$,
  \be{33.10}
	\frac{d}{dt} \io |\nabla v|^4 + c_9 \io |\nabla v|^4
	\le c_{10} \|v(\cdot,t_0)\|_{L^\infty(\Omega)}^4 \io |\nabla a|^4
	+ c_{10} \io |\nabla w|^4
	\qquad \mbox{for all } t>t_0.
  \ee
  We finally multiply the third equation in (\ref{0}) by $-\Delta w$ and use Young's inequality and (\ref{33.22}) in a straightforward
  manner to derive the inequality
  \bas
	\frac{1}{2} \frac{d}{dt} \io |\nabla w|^2
	+ \dw \io |\Delta w|^2
	+ \io |\nabla w|^2
	&=& - \io uz \Delta w \\
	&\le& \frac{\dw}{2} \io |\Delta w|^2
	+ \frac{2}{\dw} \io u^2 z^2 \\
	&\le& \frac{\dw}{2} \io |\Delta w|^2
	+ \frac{2c_1^2 c_3^2 |\Omega|}{\dw}
	\qquad \mbox{for all } t>0,
  \eas
  where a second application of (\ref{33.5}) in conjunction with (\ref{33.22}) shows that
  \bas
	\frac{\dw}{2} \io |\Delta w|^2 \ge \frac{\dw}{2c_2^2 c_6} \io |\nabla w|^4
	\qquad \mbox{for all } t>0,
  \eas
  so that, in fact,
  \be{33.11}
	\frac{d}{dt} \io |\nabla w|^2 + c_{11} \io |\nabla w|^4 + 2\io |\nabla w|^2
	\le c_{12}
	\qquad \mbox{for all } t>0
  \ee
  with $c_{11}:=\frac{\dw}{c_2^2 c_6}$ and $c_{12}:=\frac{4c_1^2 c_3^2 |\Omega|}{\dw}$.\abs
  Now in order to suitably combine (\ref{33.7}), (\ref{33.10}) and (\ref{33.11}), we abbreviate
  \bas
	b_1:=\frac{2c_8}{c_9}
	\qquad \mbox{and} \qquad
	b_2:=\frac{c_{10}b_1}{c_{11}},
  \eas
  and rely on the decay property of $v$ asserted by Lemma \ref{lem22} to fix $t_0\ge 1$ large enough such that
  \bas
	c_{10} b_1 \|v(\cdot,t_0)\|_{L^\infty(\Omega)}^4 \le c_7.
  \eas
  Then (\ref{33.7}), (\ref{33.10}) and (\ref{33.11}) imply that
  \bas
	y(t):=\io |\nabla a(\cdot,t)|^2 + b_1 \io |\nabla v(\cdot,t)|^4  + b_2 \io |\nabla w(\cdot,t)|^2,
	\qquad t\ge t_0,
  \eas
  satisfies
  \bas
	y'(t) + \frac{c_9}{2} y(t)
	&\le& \bigg\{ - c_7 \io |\nabla a|^4
	- 2 \io |\nabla a|^2
	+ c_8 \io |\nabla v|^4
	+ c_8 \bigg\} \\
	& & + b_1 \cdot \bigg\{ -c_9 \io |\nabla v|^4
	+ c_{10} \|v(\cdot,t_0)\|_{L^\infty(\Omega)}^4 \io |\nabla a|^4 + c_{10} \io |\nabla w|^4 \bigg\} \\
	& & + b_2 \cdot \bigg\{ - c_{11} \io |\nabla w|^4
	- 2 \io |\nabla w|^2
	+ c_{12} \bigg\} \\
	& & + \frac{c_9}{2} \cdot \bigg\{ \io |\nabla a|^2 + b_1\io |\nabla v|^4 + b_2\io |\nabla w|^2 \bigg\} \\[1mm]
	&=& \Big\{ - c_7 + c_{10}b_1 \|v(\cdot,t_0)\|_{L^\infty(\Omega)}^4 \Big\} \cdot \io |\nabla a|^4 \\
	& & + \Big\{ c_8-\frac{c_9}{2}b_1\Big\} \cdot \io |\nabla v|^4 \\
	& & + \big\{ c_{10}b_1-c_{11}b_2\big\} \cdot \io |\nabla w|^4 \\
	& & + \Big\{-2+\frac{c_9}{2}\Big\}\cdot \io |\nabla a|^2 \\
	& & + \Big\{-2b_2 + \frac{c_9}{2}b_2\Big\}\cdot \io |\nabla w|^2 \\
	& & + c_8 + c_{12} b_2 \\[1mm]
	&\le& c_8 + c_{12}b_2
	\qquad \mbox{for all } t>t_0,
  \eas
  because $c_9 \le 4$. As thus
  \bas
	y(t) \le c_{12}:=\max \Big\{ y(t_0) \, , \, \frac{2(c_8+c_{11}b_2)}{c_9} \Big\}
	\qquad \mbox{for all } t\ge t_0
  \eas
  by an ODE comparison argument, it particularly follows that
  \bas
	\io |\nabla v(\cdot,t)|^4 \le c_{13} :=\max \bigg\{ \sup_{s\in (0,t_0)} \io |\nabla v(\cdot,s)|^4, \frac{c_{12}}{b_1} \bigg\}
	\qquad \mbox{for all } t>0,
  \eas
  with finiteness of $c_{13}$ guaranteed by Lemma \ref{lem_basic}.
\qed
Thanks to the fact that the integrability exponent in (\ref{33.2}) exceeeds the considered spatial dimension,
through standard parabolic regularity theory this implies a uniform H\"older bound for $u$:
\begin{lem}\label{lem32}
  Suppose that $\beta>0$, $\gamma\in (0,\frac{1}{(\beta-1)_+})$, $\rho\ge 0$ and $M>0$, and that
  (\ref{init}) and (\ref{iu})-(\ref{iz}) are valid with some $\eps\in (0,\esss)$, where $\esss=\esss(\beta,\gamma,M)>0$
  is as given by Lemma \ref{lem28}.
  Then there exist $\theta\in (0,1)$ and $C>0$ such that
  \bas
	\|u\|_{C^{\theta,\frac{\theta}{2}}(\bom\times [t,t+1])} \le C
	\qquad \mbox{for all } t>0.
  \eas
\end{lem}
\proof
  We rewrite the first equation in (\ref{0}) according to $u_t=\Delta u - \nabla \cdot \psi_1(x,t) + \psi_2(x,t)$,
  $x\in\Omega, \ t>0$,
  where $\psi_1(x,t):=u(x,t)\nabla v(x,t)$ and $\psi_2(x,t):=-\rho u(x,t)z(x,t)$, $x\in\Omega, \ t>0$.
  Since (\ref{28.1}) together with the outcome of Lemma \ref{lem33} ensures that $\psi_1$ belongs to
  $L^p((0,\infty);L^q(\Omega;\R^2))$ with $p:=\infty$ and $q:=4$ satisfying $\frac{1}{p}+ \frac{2}{q}=\frac{1}{2}<1$,
  and since $\psi_2$ in bounded by (\ref{28.1}) and (\ref{28.3}),
  this directly results from well-known theory on H\"older regularity of bounded solutions to scalar parabolic equations
  (\cite{porzio_vespri}).
\qed
Straightforward interpolation between the latter and the basic convergence result from Lemma \ref{lem52} finally
yields uniform stabilization also in the first solution component:
\begin{lem}\label{lem34}
  Let $\beta>0$, $\gamma\in (0,\frac{1}{(\beta-1)_+})$, $\rho\ge 0$ and $M>0$, and suppose that
  (\ref{init}) and (\ref{iu})-(\ref{iz}) hold with some $\eps\in (0,\esss)$, where $\esss=\esss(\beta,\gamma,M)>0$
  is taken from Lemma \ref{lem28}.
  Then there exists $u_\infty>0$ such that
  \be{34.1}
	u(\cdot,t) \to u_\infty
	\quad \mbox{in } L^\infty(\Omega)
	\qquad \mbox{as } t\to\infty.
  \ee
\end{lem}
\proof
  Since Lemma \ref{lem28} together with Lemma \ref{lem25} clearly
  ensures boundedness of $(u,w,z)$ in $\Omega\times (0,\infty)$, finiteness of $\int_0^\infty \io z$ as well as
  decay to zero of $\|v(\cdot,t)\|_{L^\infty(\Omega)}$ as $t\to\infty$,
  we may invoke Lemma \ref{lem52} to find $u_\infty>0$ such that
  \be{34.11}
	u(\cdot,t)\to u_\infty
	\quad \mbox{in } L^2(\Omega)
	\qquad \mbox{as } t\to\infty,
  \ee
  and in line with Lemma \ref{lem32}, we can thereupon pick $\theta\in (0,1)$ and $c_1>0$ such that
  \be{34.2}
	\|u(\cdot,t)-u_\infty\|_{C^\theta(\bom)} \le c_1
	\qquad \mbox{for all } t>0.
  \ee
  Given $\eta>0$ we next use the compactness of the first among the continuous embeddings
  $C^\theta(\bom) \hra L^\infty(\Omega) \hra L^2(\Omega)$ to pick $c_2(\eta)>0$ such that in accordance with an associated Ehrling lemma
  we have
  \be{34.3}
	\|\varphi\|_{L^\infty(\Omega)} \le \frac{\eta}{2c_1} \|\varphi\|_{C^\theta(\bom)}
	+ c_2(\eta) \|\varphi\|_{L^2(\Omega)}
	\qquad \mbox{for all } \varphi\in C^\theta(\bom),
  \ee
  and rely on (\ref{34.11}) in verifying that for any such $\eta$ we can find $t_0(\eta)>0$ fulfilling
  \bas
	\|u(\cdot,t)-u_\infty\|_{L^2(\Omega)} \le \frac{\eta}{2c_2(\eta)}
	\qquad \mbox{for all } t>t_0(\eta).
  \eas
  Combining this with (\ref{34.3}) and (\ref{34.2}) shows that
  \bas
	\|u(\cdot,t)-u_\infty\|_{L^\infty(\Omega)}
	&\le& \frac{\eta}{2c_1}\|u(\cdot,t)-u_\infty\|_{C^\theta(\bom)}
	+c_2(\eta)\|u(\cdot,t)-u_\infty\|_{L^2(\Omega)} \\
	&\le& \frac{\eta}{2c_1} \cdot c_1
	+ c_2(\eta) \cdot \frac{\eps}{2c_2(\eta)} =\eta
	\qquad \mbox{for all } t>t_0(\eta)
  \eas
  and hence establishes (\ref{34.1}), for $\eta>0$ was arbitrary.
\qed
Accomplishing our main results now reduces to merely extracting the respectively relevant pieces of information from
the above statements:\abs
\proofc of Theorem \ref{theo23}. \quad
  Applying Lemma \ref{lem25} and Lemma \ref{lem28} to any fixed $\eps\in (0,\min\{\es,\esss\})$, with $\es=\es(\beta,\gamma)>0$ and
  $\esss=\esss(\beta,\gamma,M)>0$ as introduced there, assuming (\ref{init}), and (\ref{iu})-(\ref{iz}) we immediately obtain
  (\ref{23.4}) from Lemma \ref{lem25} and (\ref{23.3}) from Lemma \ref{lem34}, whereas
  (\ref{23.5}) and (\ref{23.6}) are direct consequences of Lemma \ref{lem28}
  when combined with Corollary \ref{cor288}.
  In the borderline case when $\rho=0$, finally, the identity $u_\infty=\ouz$ readily results from (\ref{23.4}) and the evident fact
  that $\io u(\cdot,t)=\io u_0$ for all $t>0$ by (\ref{0}).
\qed

\vspace*{10mm}
{\bf Acknowledgement.} \quad
  Youshan Tao was supported by the {\em National Natural Science Foundation of China
  (No. 11861131003)}. Michael Winkler acknowledges support of the {\em Deutsche Forschungsgemeinschaft}
  in the context of the project
  {\em Emergence of structures and advantages in cross-diffusion systems} (No.~411007140, GZ: WI 3707/5-1).

\begin{thebibliography}{99}
%
\bibitem{alemany}
  \sc Alemany, R.:
  \it Viruses in cancer treatment.
  \rm Clin. Transl. Oncol. {\bf 15}, 182-188 (2013)
\bibitem{eftimie}
  \sc Alzahrani, T., Eftimie, R., Trucu, D.:
  \it Multiscale modelling of cancer response to oncolytic viral
  therapy.
  \rm Math. Biosci. {\bf 310}, 76-95 (2019)
\bibitem{bischoff}
  \sc Bischoff, J.R., Kirn, D.H., Williams, A., Heise, C., Horn, S.,  Muna, M.,  Ng, L., Nye, J.A.,
      Sampson-Johannes, A., Fattaey, A., McCormick, F.:
  \it An adenovirus mutant that replicates selectively in p53-deficient
      human tumor cells.
  \rm Science {\bf 274}, 373-376 (1996)
\bibitem{cao_zamp2016}
  \sc Cao, X.:
  \it Boundedness in a three-dimensional chemotaxis-haptotaxis system.
  \rm Z.~ Angew.~ Math.~ Phys. {\bf 67}, 11 (2016)
\bibitem{zhen_chen}
  \sc Chen, Z.:
  \it Dampening effect of logistic source in a two-dimensional haptotaxis system with nonlinear zero-order interaction.
  \rm Preprint
\bibitem{coffey}
  \sc Coffey, M.C., Strong, J.E., Forsyth, P.A., Lee, P.W.K.:
  \it Reovirus therapy of tumors with activated Ras pathways.
  \rm Science {\bf 282}, 1332-1334 (1998)
\bibitem{fontelos_friedman_hu}
  \sc Fontelos, M.A., Friedman, A., Hu, B.:
  \it Mathematical analysis of a model for the initiation of angiogenesis.
  \rm SIAM J. Math. Anal. {\bf 33}, 1330-1355 (2002)
\bibitem{friedman_tello}
  \sc Friedman, A., Tello, J.I.:
  \it Stability of solutions of chemotaxis equations in reinforced random walks.
  \rm J. Math. Anal. Appl. {\bf 272}, 138-163 (2002)
\bibitem{hillen_painter_win_M3AS}
  \sc Hillen, T., Painter, K.J., Winkler, M.:
  \it Convergence of a cancer invasion model to a logistic chemotaxis model.
  \rm Math. Mod. Meth. Appl. Sci.  {\bf 23}, 165-198 (2013)
\bibitem{jain}
  \sc Jain, R.:
  \it Barriers to drug delivery in solid tumors.
  \rm Sci. Am. {\bf 271}, 58-65 (1994)
\bibitem{komarova}
  \sc Komarova, N.L.:
  \it  Viral reproductive strategies: how can lytic viruses be evolutionarily competitive?
 \rm  J.~Theor.~Biol. {\bf 249}, 766-784 (2007)
\bibitem{yifu_wang}
  \sc Li, J., Wang, Y.:
  \it Boundedness in a haptotactic cross-diffusion system modeling oncolytic virotherapy.
  \rm Preprint
\bibitem{litcanu_cmr}
  \sc Li\c{t}canu, G., Morales-Rodrigo, C.:
  \it Asymptotic behavior of global solutions to a model of cell invasion.
  \rm Math. Models Methods Appl. Sci. {\bf 20}, 1721-1758 (2010)
\bibitem{martuza}
  \sc Martuza, R.L., Malick, A., Markert, J.M., Ruffner, K.L., Coen, D.M.:
  \it Experimental therapy of human glioma by means of a genetically engineered virus mutant.
  \rm Science {\bf 252}, 854-856 (1991)
\bibitem{morales_tello}
  \sc Morales-Rodrigo, C., Tello, J.I.:
  \it Global existence and asymptotic behavior of a tumor angiogenesis model with chemotaxis and haptotaxis.
  \rm Math. Models Methods Appl. Sci. {\bf 24}, 427-464 (2014)
\bibitem{painter_maini_othmer}
  \sc Painter, K.J., Maini, P.K., Othmer, H.G.:
  \it Stripe formation in juvenile \textit{Pomacanthus} explained by a
      generalized Turing mechanism with chemotaxis.
  \rm Proc. Natl. Acad. Sci. USA {\bf 96}, 5549-5554 (1999)
\bibitem{pang_wang_m3as2018}
  \sc Pang, P.Y.H., Wang, Y.:
  \it Global boundedness of solutions to a chemotaxis-haptotaxis model with tissue remodeling.
  \rm Math.~Mod.~Meth.~Appl.~Sci. {\bf 28}, 2211-2235 (2018)
\bibitem{porzio_vespri}
 \sc Porzio, M.M., Vespri, V.:
  \it H\"{o}lder estimates for local solutions of some doubly nonlinear degenerate parabolic equations.
  \rm J.~Differential Equations {\bf 103} (1), 146-178 (1993)
\bibitem{pruess_mmnp}
  \sc Pr\"uss, J., Zacher, R., Schnaubelt, R.:
  \it Global asymptotic stability of equilibria in models for virus dynamics.
  \rm Math.~Model.~Nat.~Phenom. {\bf 3} (7), 126-142 (2008)
\bibitem{rodriguez_win}
  \sc Rodriguez, N., Winkler, M.:
  \it On the global existence and qualitative behavior of
  one-dimensional solutions to a model for urban crime.
  \rm Preprint
\bibitem{russell}
  \sc Russell, S.J., Peng, K.-W., Bell, J. C.: \it Oncolytic
  virotherapy.
  \rm Nature Biotechnology {\bf 30}, 658-670 (2012)
\bibitem{stinner_surulescu_winkler_SIMA}
  \sc Stinner, C., Surulescu, C., Winkler, M.:
  \it Global weak solutions in a PDE-ODE system modeling multiscale cancer cell invasion.
  \rm SIAM J. Math. Anal. \textbf{46},  1969-2007 (2014)
\bibitem{surulescu_win_CMS}
  \sc Winkler, M., Surulescu, C.:
  \it A global weak solutions to a strongly degenerate haptotaxis model.
  \rm Commun. Math. Sci. {\bf 15}, 1581-1616 (2017)
\bibitem{swabb}
  \sc Swabb, E.A., Wei, J., Gullino, P.M.:
  \it Diffusion and convection in normal and neoplastic tissues.
  \rm Cancer Res. {\bf 34}, 2814-2822 (1974)
\bibitem{xueyan_tao}
  \sc Tao, X.: \it Global classical solutions to an oncolytic viral therapy model with triply haptotactic terms in 1D.
  \rm Preprint
\bibitem{tao_wang_sima2009}
  \sc Tao, Y., Wang, M.:
  \it A combined chemotaxis-haptotaxis system: The role of logistic
  source.
  \rm   SIAM J. Math. Anal. {\bf 41}, 1533-1558  (2009)
\bibitem{taowin_JDE2014}
  \sc Tao, Y., Winkler, M.:
  \it Energy-type estimates and global solvability in a two-dimensional chemotaxis-haptotaxis
  model with remodeling of non-diffusible attractant.
  \rm J.~Differential Eq. {\bf 257}, 784-815 (2014)
\bibitem{taowin_NON}
  \sc Tao, Y., Winkler, M.:
  \it Dominance of chemotaxis in a chemotaxis-haptotaxis model.
  \rm Nonlinearity {\bf 27} (6), 1225-1239 (2014)
\bibitem{taowin_sima2015}
  \sc Tao, Y., Winkler, M.:
  \it Large time behavior in a mutidimensional chemotaxis-haptotaxis model with slow signal diffusion.
  \rm SIAM J.~ Math.~ Anal. {\bf 47}, 4229-4250 (2015)
\bibitem{taowin_261}
  \sc Tao, Y., Winkler, M.:
  \it Global classical solutions to a doubly haptotactic cross-diffusion system modeling oncolytic
      virotherapy.
  \rm J. Differential Equations {\bf 268}, 4973-4997 (2020)
\bibitem{taowin_262}
  \sc Tao, Y., Winkler, M.:
  \it Critical mass for infinite-time blow-up in a haptotaxis system with nonlinear zero-order
  interaction.
  \rm Discr.~Cont.~Dyn.~Syst.~A, doi: 10.3934/dcds.2020216
\bibitem{taowin_263}
  \sc Tao, Y., Winkler, M.:
  \it A critical virus production rate for blow-up suppression in a haptotatxis model for oncolytic
      virotherapy.
  \rm Nonlinear Analysis {\bf 198}, Art. 111870 (2020)
\bibitem{taowin_265}
  \sc Tao, Y., Winkler, M.:
  \it A critical virus production rate for efficiency of oncolytic virotherapy
  \rm Eur. J. Appl. Math., doi:10.1017/S0956792520000133
\bibitem{taowin_crime}
  \sc Tao, Y., Winkler, M.:
  \it Global smooth solutions in a two-dimensional cross-diffusion system modeling propagation of urban crime.
  \rm Preprint
\bibitem{vaha}
  \sc V\"{a}h\"{a}-Koskela, M., Hinkkanen, A.:
  \it Tumor restrictions to oncolytic virus.
  \rm Biomedicines {\bf 2} (2), 163-194 (2014)
\bibitem{walker_webb}
  \sc  Walker, C., Webb, G.F.:
  \it Global existence of classical solutions for a haptotaxis model.
  \rm SIAM J. Math. Anal. {\bf 38}, 1694-1713 (2007)
\bibitem{win_jmpa2018}
  \sc Winkler, M.:
  \it Singular structure formation in a degenerate haptotaxis model involving myopic diffusion.
  \rm J. Math. Pures Appl. {\bf 112}, 118-169 (2018)
\bibitem{wong}
  \sc Wong, H., Lemoine, N., Wang, Y.:
  \it Oncolytic viruses for cancer therapy: overcoming the obstacles.
  \rm Viruses {\bf 2} (1), 78-106 (2010)
\bibitem{zhigun_surulescu_uatay}
  \sc Zhigun, A., Surulescu, C., Uatay, A.:
  \it Global existence for a degenerate haptotaxis model of cancer invasion.
  \rm Z. Angew. Math. Phys. {\bf 67},  Art. 146, 29 pp. (2016)
%
\end{thebibliography}
\end{document}